\documentclass[12pt,draft]{extartic}

\usepackage{amsmath,amsthm,amssymb,calc,mathrsfs,amsfonts}
\usepackage[dvips]{graphics, graphicx}

\usepackage[english]{babel}

\def\eps{\varepsilon}
\def\gs{\geqslant}
\def\ls{\leqslant}

\newcounter{num}[section]

\newcommand{\Th}{\refstepcounter{num}
{\bf Theorem \arabic{section}.\arabic{num} }}

\newcommand{\Lemma}{\refstepcounter{num}
{\bf Lemma \arabic{section}.\arabic{num} }}

\newcommand{\con}{\refstepcounter{num}
{\bf Conjecture \arabic{section}.\arabic{num} }}

\newcommand{\Cor}{\refstepcounter{num}
{\bf Corollary \arabic{section}.\arabic{num} }}

\newcommand{\Def}{\refstepcounter{num}
{\it Definition \arabic{section}.\arabic{num} }}

\newcommand{\Proof}{{\bf Proof. }}

\def\eps{\varepsilon}
\def\Z{\mathbf{Z}}

\def\_phi{\varphi}
\def\d{\delta}

\def\eps{\varepsilon}
\def\_phi{\varphi}
\def\k{\kappa}
\def\e{\eta}
\def\a{\alpha}
\def\d{\delta}

\def\F{\widehat}
\def\L{\Lambda}

\def\t{\tilde}

\def\ov{\overline}
\def\z{{\mathbb Z}}
\def\C{{\mathbb C}}

\def\T{{\mathbb T}}

\def\Z_N{{\mathbb Z}_N}
\def\Z{{\mathbb Z}}

\def\T{\mathcal T}
\def\s{\bf g}

\def\f{{\mathbb F}}
\def\Gr{{\mathbf G}}

\def\dim{{\rm dim }}

\def\G{\Gamma}
\def\g{\gamma}

\def\l{\left}
\def\r{\right}

\def\Spec{{\rm Spec\,}}

\def\sbeq{\subseteq}
\newcommand{\lf}{\left\lfloor}
\newcommand{\rf}{\right\rfloor}

\def\no{\noindent}

\author{\sc Tomasz Schoen\footnote{The author is partly supported by MNSW grant N N201
543538.} ~ and ~ Ilya~D.~Shkredov\footnote{The author is supported
by grant RFFI NN 06-01-00383, 11-01-00759 and grant Leading Scientific Schools No. 8684.2010.1}}
\title{\bf Roth's theorem in  many variables}

\date{}
\begin{document}
\maketitle

\begin{center}
    Abstract
\end{center}

{\it \small
We prove, in particular, that if   $A\sbeq \{1,\dots ,N\}$ has no nontrivial solution to the equation
    $x_1+x_2+x_3+x_4+x_5=5y$ then  $|A|\ll Ne^{-c(\log N)^{1/6-\eps}},\, c>0$.
    In view of the well-known Behrend construction this estimate is close to  best possible.
}
\\

\section{Introduction}
\label{sec:introduction}

 The celebrated theorem of Roth \cite{Roth_AP3} asserts that every subset of $\{1,\dots, N\}$ that does
 not contain any  three term arithmetic progression has size  $O(N/\log\log N)$.
 There are numerous refinements of  Roth's result \cite{Bu,Bourgain_AP2007,a_H-B,Sz_inf}. Currently
the best known upper bound $O( N/(\log N)^{1-o(1)})$ is due to
Sanders \cite{Sanders_Roth_1-eps}.
The comprehensive history of the subject can be found in \cite{Tao_Vu_book}.

It turns out that the Roth's method gives a similar upper bound for the size of sets
having no nontrivial solutions to a
invariant
 linear equation
\begin{equation}\label{e}
a_1x_1+\dots+a_kx_k=0,
\end{equation}
i.e.  $a_1+\dots+a_k=0,\, k\gs 3$ (three term arithmetic progressions corresponds to the equation $x+y=2z$).
On the other hand, the well-known construction of Behrend \cite{Behrend,Elkin,Koester,Moser,Ra,Salem1,Salem2}
provides large sets having no solution to  certain kind of invariant equations.
He showed that there are
subsets of $\{1,\dots,N\}$ of size $Ne^{(-C_{b,k}\sqrt{\log N})}$ without solution to
the invariant equation
\begin{equation}\label{eq: Behrend}
a_1x_1+\dots +a_kx_k=by,
\end{equation}
where  $a_1+\dots+a_k=b,~a_i>0.$

The aim of this paper is to establish a new upper bound for subsets of $\{1,\dots,N\}$ having no solution
to an invariant equation in at least $6$ variables.

\bigskip

\Th
{\it
    Let $N$ and  $k\gs 6$ be positive  integers.
    Let $A \subseteq \{1,\dots,N\}$ be a set having no solution to
     the equation (\ref{e}), where all $x_1,\dots,x_k,y$ are distinct integers.
    Then
    \begin{equation}\label{est:cardinality_A_5y''_intr}
        |A|\ll {\exp\Big (-c\Big (\frac{ \log N}{\log \log N}\Big )^{1/6} \Big )}N \,,
    \end{equation}
where $c=c(a_1,\dots,a_k).$}
\label{t:Roth_Schoen_Gr''_intr}

\bigskip

Observe, that Theorem \ref{t:Roth_Schoen_Gr''_intr} together with Behrend's example give
a reasonable  estimates for all equations of the type (\ref{eq: Behrend}).
Let us also formulate an immediate corollary to Theorem \ref{t:Roth_Schoen_Gr''_intr}
for the equation
\begin{equation}\label{f}
x_1+x_2+x_3+x_4+x_5=5y
\end{equation}
which is very close to the most intriguing  case
$x+y=2z$.

\bigskip

\Cor
{\it Suppose that $A\sbeq \{1,\dots,N\}$ has no solution to the equation (\ref{f}) with distinct integers.
 Then there exists a constant $c>0$ such that
\begin{equation*}
        |A|\ll {\exp\Big (-c\Big (\frac{ \log N}{\log \log N}\Big )^{1/6} \Big)}N \,.
\end{equation*}
}
\label{c:roth5x}

Our argument heavily relies on a recent work on  Polynomial Freiman-Ruzsa Conjecture, by Sanders
\cite{Sanders_2A-2A_new_optimal} (see also  \cite{sch}).
A fundamental tool in our approach is a
version of Bogolyubov-Ruzsa Lemma proved in \cite{Sanders_2A-2A_new_optimal}.
 We also use  the density
  increment method introduced by Roth, however in a different way.
The density increment is not deduced  from the existence of a large Fourier coefficient
of a set $A,\, |A|=\a N,$ having no solution to an equation (\ref{e}) (which is always the case).
We will be rather interested in finding
a translation of a  large Bohr set in $a_1\cdot A+a_2\cdot A+a_3\cdot A+a_4\cdot A,$ from which
one can easily deduce a large density increment of $A$ by a constant factor on some large Bohr set.
By Sanders' theorem dimension of the  Bohr set increases by $O(\log^4(1/\a))$ in each iteration step, which makes the
argument very effective.

The paper is organized as follows. We start  with proving analogues of Theorem
\ref{t:Roth_Schoen_Gr''_intr} and Corollary \ref{c:roth5x} for finite fields in section \ref{sec:F_p^n}.
The argument is especially simple and transparent in this case.
 Theorem \ref{t:Roth_Schoen_Gr''_intr} is proved in next three sections.
In section \ref{sec:Bohr_sets} we recall some basic properties of Bohr sets in abelian groups.
In section \ref{sec:Sanders_Bohr} we prove a local version of Sanders result.
The next section contains the proof of  Theorem \ref{t:Roth_Schoen_Gr''_intr}.
We conclude the paper with a discussion concerning consequences of
Polynomial Freiman--Ruzsa Conjecture
for sets having no solutions to an invariant linear equation with distinct integers.

\section{Notation}
\label{sec:notation}

Let $\Gr=(\Gr,+)$ be a finite Abelian group with additive group
operation $+$, and let $N=|\Gr|$.
 By $\F{\Gr}$ we denote the
Pontryagin dual of $\Gr$, i.e. the space of
homomorphisms $\gamma$ from $\Gr$ to $S^1$. It is well known that $\F{\Gr}$ is an
additive group which is isomorphic to $\Gr$. The
Fourier coefficients of $f: \Gr\rightarrow \C$ are defined by
$$\F{f}(\g)=\sum_{x\in \Gr}f(x)\ov{\g}(x).$$
By the  convolution of two function $f,g:\Gr\rightarrow \C$ we mean
$$(f*g)(x)=\sum_{x\in \Gr}f(x)g(y-x).$$ It is easy to see that
$\F{f*g}(\g)=\F{f}(\g) \F{g}(\g).$ If $X$ is a nonempty set, then by $\mu_X$ we denote the uniform probability measure on $X$ and  let
$$ \Spec_\epsilon (\mu_X) := \{ \gamma \in \F{\Gr} ~:~ |\F{X} (\gamma)| \gs \epsilon |X| \}.$$

Let
$\Z_p = \Z/p\Z$, and
$\f_p^* = \Z_p \setminus \{ 0 \}$.
If $A$ is a set, then we write $A(x)$ for its characteristic function i.e. $A(x) = 1$ if $x\in A$ and $A(x)=0$ otherwise.
All logarithms  are to  base $2.$
The signs $\ll$ and $\gg$ are  usual Vinogradov's symbols.

\section{Finite fields model}
\label{sec:F_p^n}


\bigskip

In this section we present proofs of Corollary \ref{c:roth5x} and Theorem \ref{t:Roth_Schoen_Gr''_intr}
  in finite fields setting. Here we assume that  $a_1,\dots,a_k\in \f^*_p.$  The case of  $\f^n_p,$
in view of its linear  space structure over $\f_p$, is considerable simpler than the case of $\Z.$
Even the simplest version of Roth's argument yields an estimate
 $O_p(p^n/n^{k-2})$
for size of sets free of  solution to (\ref{e})
  (see \cite{Meshulam,LS09}, \cite{Sanders_Z_4^n},\cite{Sanders_log}).

Our main tool is the following finite fields version of Sanders' theorem \cite{Sanders_2A-2A_new_optimal}.

\bigskip

\Lemma
{\em
Suppose that $A,S\subseteq \f_p^n$ are finite non--empty sets such that
$|A+S| \le K \min \{ |A|, |S| \}$.
Then $A-A+S-S$ contains a subspace $V$ of codimension  at most
$O_p(\log^{4} K)$.
}
\label{l:Sanders_F_p^n}

\bigskip

The proof of the next theorem illustrates  the main idea of our approach.

\bigskip

\Th
{\em Suppose that $A\sbeq \f_p^n,\, p\not=5,$
and $A$  has no nontrivial solution to (\ref{f}) with $x_i\not=y$ for some $i$.
Then
$$
    |A|\ls p^n \cdot \exp(-c_p(\log p^n)^{1/5}) \,
$$
for some positive constant $c_p.$}
\label{t:5y_F_p^n}

\proof
Suppose that $A\sbeq \f_p^n$ has density $\a$ and contains
no solution to $(\ref{f}).$
We split $A$ into two disjoint sets
$A_1,A_2$ of equal size. Clearly, there exists  $z\in \f^n_p$
such that
$$|A_1\cap (z-A_2)|\gg  \a^2p^n\,.$$
Let us put $B=A_1\cap (z-A_2)$.

By Lemma \ref{l:Sanders_F_p^n}, there exists a subspace $V$ of codimension at most
$O_p(\log^{4} (1/\a))$
such that $V\sbeq 2B-2B$, so that
$$
    2z+V \sbeq 2A_1+2A_2 \,.
$$
Therefore, in view of $A_1\cap A_2=\emptyset$,  we have
$5y-x\not\in 2z+V$ for all $x,y\in A,$ hence  $A$
intersects at most half of cosets of $V,$ which implies
$$|A\cap (v+V)|\ge 2\a|V|,$$
for some $v.$ Thus, $(A-v)\cap V$ is free of solutions to
(\ref{f}) and has density at least $2\a$ on $V.$ After $t$
iterations we obtain a subspace of codimension at most
$O_p(t \cdot \log^{4} (1/\a))$ such that
$$
    |(A-v_t)\cap V_t|\ge 2^t\a|V_t| \,,
$$
for some $v_t.$ Since the density is always at most one we can
iterate this procedure at most $\log(1/\a)+1$ times. Hence
$$
    (\log(1/\a)+1) \cdot \log^{4} (1/\a) \gg_p n \,,
$$
so that
$$ \a\ls \exp(-c_pn^{1/5}) \,$$
for some positive constant $c_p.\hfill\Box$
\bigskip

 To prove the main result of the section we will need the following  consequence of Lemma \ref{l:Sanders_F_p^n}.
We sketch  its proof here; the interested reader will find  all details in Section \ref{sec:Sanders_Bohr}.

\bigskip

\Lemma
{\em
Let $A_1,\dots, A_k\sbeq \f^n_p$ be sets of density at least $\a$.
Then
$A_1-A_1+\dots +A_k-A_k$
contains  a  subspace $V$ of codimension  at most
$O_p(k^{-3}\log^{4} (1/\a))$.}
\label{l:trick_small_sumsets}

\medskip

\proof
We have
$$|A_1|\le |A_1+A_2|\le\dots \le |A_1+\dots +A_k|\le \a^{-1}|A_1|, $$
so that there exists $2\le i\le k$ such that
$$|A_1+\dots +A_i|\le \a^{-1/(k-1)}|A_1+\dots +A_{i-1}|.$$
Thus, setting  $A=A_1+\dots+A_{i-1},\, S=T=A_i$, we have  $|A+S|\ls  \a^{-1/(k-1)}|A|$,
$|S+T|\ls \a^{-1}|S|.$ Then applying Lemma \ref{l:Sanders_F_p^n} and Theorem \ref{t:Sanders_2A-2A_reformulation} (see Section \ref{sec:Sanders_Bohr})
  we infer that there is a subspace $V$ of codimension
$O_p(\log^3(\a^{-1/(k-1)})\cdot \log (1/\a))=O_p(k^{-3}\log^{4} (1/\a))$
such that
$$
    v+V\sbeq A_1-A_1+\dots +A_i-A_i \sbeq A_1-A_1+\dots +A_k-A_k \,,
$$
and the assertion follows. $\hfill\Box$

\bigskip

\Th
{\em Suppose that $A\sbeq
\f_p^n$ has no solution   with distinct elements to an invariant
equation
\begin{equation}\label{equation:F_p^n}
a_1x_1+\dots+a_kx_k=0,
\end{equation}
where  $a_1,\dots,a_k\in \f_p^*$ and $k\gs 6.$
Then
$$
    |A|\ls k p^n \cdot \exp(-c_p(k^3\log p^n)^{1/5}) \,.
$$
for a positive constant $c_p$.}
\label{t:5y_general_F_p^n}

\medskip

\proof
Suppose  $A\sbeq \f_p^n$ has no solution
 with distinct elements to (\ref{equation:F_p^n})
and $|A|=\a p^n.$
Let  $A_1,\dots, A_{2l},\, l=\lf (k-2)/2 \rf$ be arbitrary disjoint
subsets of $A$ of size  $\lf |A|/(5k)\rf$ and put $A'=A\setminus
\bigcup A_i.$
Clearly,  there are $z_1,\dots,z_{l}$ such that
$$
    |(a_{2i-1}\cdot A_{2i-1})\cap(z_i-a_{2i}\cdot A_{2i})|\gg (k/\a)^{2}p^n
    $$
 and let $B_i,\, 1\ls i\ls l,$   be the sets on the left hand
side in the above inequalities, respectively.
By Lemma \ref{l:trick_small_sumsets}, applied for $B_1,\dots, B_l$
and $K=O((k/\a)^2)$ there is a subspace $V$ of codimension
$d=O_p(k^{-3}\log^{4} (k/{\a}))$   such
that
$$
    V\sbeq B_1-B_1 +\dots+B_l-B_l \,,
$$
so that
$$v+V\sbeq a_1\cdot A_1+\dots+a_{k-2}\cdot A_{k-2}$$
for some $v.$ Since $A$ does not contain any solution to (\ref{equation:F_p^n})
with distinct elements it follows that
$$a_{k-1}x+a_ky\notin v+V,$$
for all $x,y\in  A', \, x \not=y.$
Hence, if for some $w$  the coset $w+V$ contains at least $2$
elements of $A',$ then $-a^{-1}_k(a_{k-1}w-v)+V$ is disjoint from
$A'$. The number of cosets of $V$ sharing exactly $1$ element with
$A$ is trivially at most $ p^d.$ Thus, there exists $w'$ such that
$|A'\cap (w'+V)|\gs  \frac{(4/5)\a p^n-p^d}{p^d/2} |V|,$ which is
at least $(3/2)\a |V|,$ provided that
\begin{equation}\label{d}
    p^{n-d}\gg \a^{-1}.
\end{equation}
After $t$ iterates of this argument we obtain a subspace $V_t$ of
codimension  $O_p(t k^{-3}\log^{4} (k/{\a}))$ such that
$$|(A-v_t)\cap V_t|\gs (3/2)^t\a|V_t|.$$
Since $(3/2)^t\a\ls 1$ it follows that  $t \ls 2\log (1/\a)$.
Thus, (\ref{d}) must be violated after at most $2\log(1/\a)$
steps, in particular $p^{n-2\log(1/\a)d} \ll \a^{-1},$ so that
$$
    k^{-3}\log(1/\a)  \log^{4} ({k}/{\a}) \gg_p n/2 \,.
$$
Hence
$\a \ls k \exp(-c_p(k^3\log p^n)^{1/5}) \,.\hfill \Box$


\section{Basic properties of Bohr sets}
\label{sec:Bohr_sets}

\bigskip

Bohr sets were introduced to additive number theory by Ruzsa
\cite{Ruzsa_freiman}. Bourgain \cite{Bu} was the first, who used
Fourier analysis on Bohr sets to improve estimate in Roth's theorem.
Sanders \cite{Sanders_2A-2A_new_optimal} further developed the
theory of Bohr sets proving many important theorems, see for example
Lemma \ref{l:Chang+large_Bohr_exp_sums} below.

 Let $\G$ be a subset of $\F{\Gr}$, $|\G| = d$,
 and $\eps = (\eps_1,\dots,\eps_d) \in (0,1]^d$.

\bigskip

\Def Define the Bohr set  $B = B(\G, \eps)$ setting
\[
  B(\G, \eps) = \{ n \in \Gr ~:~ \| \g_j(n)\| < \eps_{j} \mbox{ for all } \g_j \in \G \} \,,
 \] where $\|x\|=|\arg x|/2\pi.$
\bigskip

The number  $d$~ is called  {\it dimension } of  $B$ and is denoted
by $\dim B$. If $M = B + n$, $n\in \Gr,$ is a translation of a Bohr
set  $B$, we put $\dim M = \dim B$. The {\it intersection} $B\wedge
B'$ of two Bohr sets $B = B(\G,\eps)$ and $B' = B(\G',\eps')$ is the
Bohr set with the generating set $\G\cup \G'$ and new vector
$\t{\eps}=\min (\eps_j,\eps'_j)$. We write ${B}' \ls {B}$ for two
Bohr sets ${B} = B({\G},{\eps})$, ${B}' = B({\G}',{\eps}')$ if ${\G}
\subseteq {\G}'$ and ${\eps}'_j \ls {\eps}_j$, $j\in [\dim B]$. Thus
${B}' \le {B}$ implies that ${B}' \subseteq {B}$ and always $B\wedge
B' \le B,B'$. Furthermore, if $B=B(\G,\eps)$ and $\rho>0$ then by
$B_\rho$ we mean $B(\G,\rho \eps).$

\bigskip

\Def  A Bohr set $B = B (\G, \eps)$ is called {\it regular},
if for every  $\eta,\, d|\eta|\ls  1/100$
we have
\begin{equation}\label{cond:reg_size}
  (1-100d|\eta|)|B_1| < { |B_{1+\eta}| } < (1+100d|\eta|)|B_1| \,.
\end{equation}

\bigskip

We  formulate  a sequence of basic properties of Bohr (see \cite{Bu}), which will be used later.

\bigskip

\Lemma {\it Let  $B (\G,\eps)$ be
a Bohr set. Then there exists $\eps_1$ such that $
  \frac{\eps}{2} < \eps_1 < \eps
$ and   $B (\G,\eps_1)$ is regular.
\label{l:Reg_B}
}

\bigskip

\Lemma
{\it
Let $B (\G,\eps)$ be a Bohr set. Then
\[
  |B (\G,\eps)| \ge \frac{N}{2} \prod_{j=1}^d \eps_j \,.
\]
}
\label{l:Bohr_est}

\bigskip

\Lemma
{\it
    Let $B (\G,\eps)$ be a Bohr set.
    Then
    $$
        |B(\G,\eps)| \ls 8^{|\G|+1} |B(\G,\eps/2)| \,.
    $$

}
\label{l:entropy_Bohr}

\bigskip
\Lemma
{\it
    Suppose that $B^{(1)}, \dots, B^{(k)}$ is a sequence of Bohr sets.
    Then
    $$
        \mu_{\Gr} (\bigwedge_{i=1}^k B^{(i)}) \ge \prod_{i=1}^k \mu_{\Gr} (B^{(i)}_{1/2}) \,.
    $$
}
\label{l:Bohr_intersection_Sanders}

The next lemma is  due to  Bourgain \cite{Bu}. It shows the
fundamental property of regular Bohr sets. We recall  his argument
for the sake of completeness.
\bigskip

\Lemma {\it Let
$B = B (\G, \eps)$ be a regular Bohr set. Then for every Bohr set $B' = B (\G,
\eps')$ such that $\eps'\ls \k \eps/ (100 d)$ we have: \\
$1)~$ the number of $n's$ such that $( B * B' ) (n) > 0$ does not exceed  $|B|(1+\kappa)$,\\
$2)~$ the number of $n's$ such that $( B * B') (n) = |B'|$ is
greater than $|B|(1-\kappa)$ and
\begin{equation}\label{2k}
  \l\| (\mu_B * \mu_{B'})(n) - \mu_B (n) \r\|_1 < 2\kappa \,.
\end{equation}
\label{l:L_pm} } \proof If $(B * B^{'}) (n) > 0,$ then there exists
$m$ such that for any $\g_j \in \G$, we have
\begin{equation*}\label{tM2}
    \| \g_j \cdot m \| < \frac{\k}{100d} \eps_j , \quad  \| \g_j \cdot (n-m) \| < \eps_j \,,
\end{equation*}
so that
\begin{equation*}\label{}
  \| \g_j \cdot n \|
             < \Big( 1 + \frac{\k}{100d} \Big) \eps_j \,,
\end{equation*}
for all $\g_j \in \G$. Therefore
$
  n\in B^{+} := B \left( \G, \left( 1 + \frac{\k}{100d} \right) \eps \right)  \,
$
and by Lemma \ref{l:Reg_B}, we have $|B^{+}| \ls (1+\k) |B|$.

On the other hand, if
\begin{equation*}\label{L-}
  n\in B^{-} := B \left( \G, \left( 1 - \frac{\k}{100d} \right) \eps
  \right)\,,
\end{equation*}
then $(B * B^{'}) (n) = |B'|$. Using Lemma \ref{l:Reg_B}, we
obtain $|B^{-}| \ge (1-\k) |B|$.

To prove (\ref{2k}) observe that
$$
    \l\| (\mu_B * \mu_{B'})(n) - \mu_B (n) \r\|_1
    =
    \l\| (\mu_B * \mu_{B'})(n) - \mu_B (n) \r\|_{ l^1 ( B^{+} \setminus B^{-}) }
  \ls
  \frac{|B^{+}| - |B^{-}|}{|B|} < 2 \k
$$
as required. $\hfill\Box$

\bigskip

\Cor {\it With the assumptions of  Lemma  \ref{l:L_pm} we have  $|B| \ls |B + B'| \ls |B^+|\ls
(1+\k) |B|$. }

\bigskip

Notice  that for every $\g\in \Z^*_p$ and a  Bohr set $B(\G,\eps)$ we have
 $\g \cdot B(\G,\eps) = B(\g^{-1}\cdot \G, \eps)$.
Thus, if $B(\G,\eps)$ is a regular,  then
 $\g \cdot B(\G,\eps)$ is regular as well.

\section{A variant of Sanders' theorem}
\label{sec:Sanders_Bohr}

Very recently Sanders \cite{Sanders_2A-2A_new_optimal}  proved the following
remarkable result.

\bigskip

\Th
{\it
    Suppose that $\Gr$ is an abelian group and  $A,S\subseteq \Gr$ are finite non--empty sets such that
    $|A+S| \ls K \min \{ |A|, |S| \}$.
    Then $(A-A)+(S-S)$ contains a proper symmetric $d(K)$--dimensional coset progression $M$
    of size $\exp (-h(K)) |A+S|$.
    Moreover, we may take $d(K) = O(\log^6 K)$, and $h(K) = O(\log^6 K \log \log K)$.
}
\label{t:Sanders_2A-2A_new_optimal}

\bigskip

The aim of this section is to show the following modification of
Sanders' theorem which is  crucial for our argument.

\bigskip

\Th
{\it
    Let $\eps,\delta \in (0,1]$ be real numbers.
    Let $A,A'$ be subsets of a regular Bohr sets $B$ and let
    $S,S'$ be subsets of a regular Bohr sets $B_\eps$, where $\eps \ls 1/(100d)$ and $d=\dim B$.
    Suppose that $\mu_B (A), \mu_{B} (A'), \mu_{B_\eps} (S), \mu_{B_\eps} (S') \gs \a$.
    Then the set $(A-A')+(S-S')$ contains a translation of a regular Bohr set $z+\t{B}$ such that
     $\dim \t{B} = d + O(\log^4 (1/\a))$ and
    \begin{equation}\label{f:t(B)_size}
        |\t{B}| \gs \exp (-O(d\log d + d\log (1/\eps) + \log^4 (1/\a)\log d +\log^5 (1/\a)+d\log(1/\a))) |B|
        \,.
    \end{equation}

}
\label{t:Sanders_2A-2A_reformulation}

\bigskip

Observe that the statement above with $O(d^4 + \log^4 (1/\a))$ instead
of $d + O(\log^4 (1/\a))$ is a direct consequence of Theorem
\ref{t:Sanders_2A-2A_new_optimal} (see the beginning of the proof
of Theorem \ref{t:Sanders_2A-2A_reformulation}).

Next we will  formulate two results, which will be used in the
course of the proof of Theorem \ref{t:Sanders_2A-2A_reformulation}.
The first lemma, proved by Sanders \cite{Sanders_2A-2A_new_optimal},
is a version of Croot--Sisask theorem
\cite{Croot_Sisask_convolutions}.

\bigskip

\Lemma
 {\it
    Suppose that $\Gr$ is a group, $A,S,T \subseteq \Gr$ are finite non--empty sets such that
    $|AS| \ls K|A|$ and $|TS|\ls L|S|$.
    Let $\epsilon \in (0,1]$ and let  $h$ be a positive integer.
    Then there is $t\in T$ and a
    set
    $X\subseteq T-t$,
    with
    $$
        |X| \gs \exp (-O(\epsilon^{-2} h^2 \log K \log L)) |T|
    $$
    such that
    $$
        | \mu_{A^{-1}} * AS * \mu_{S^{-1}}(x) - 1| \ls \epsilon \quad \mbox{ for all } \quad x \in X^h \,.
    $$
}
\label{pr:Periodic_convolution_3_sets}

\bigskip

The next lemma is a special case of Lemma 5.3 from
\cite{Sanders_2A-2A_new_optimal}. This is a local version of Chang's
spectral lemma \cite{Chang}, which is  another important result
recently proved in additive combinatorics.
\bigskip

\Lemma
{\it
    Let $\epsilon, \nu, \rho $ be positive real number.
    Suppose that $B$ is a regular Bohr set and let $X\subseteq B.$
    Then there is a set $\L$ of size $O(\epsilon^{-2} \log (2\mu^{-1/2}_B (X)))$
    such that for any $\gamma\in \Spec_\epsilon (\mu_X)$
    we have
    $$
        |1-\gamma(x)| = O(|\L| (\nu + \rho \dim^2 (B)) ) \quad \mbox{ for all} \quad x\in B_{\rho} \wedge B'_\nu \,,
    $$
    where $B' = B(\L,1/2)$.
}
\label{l:Chang+large_Bohr_exp_sums}

\bigskip

{\it Proof of Theorem \ref{t:Sanders_2A-2A_reformulation}} Applying
Lemma \ref{pr:Periodic_convolution_3_sets} with $A,\, S$ and  $T=B_{\d},\,
\d={\eps/100d}$ and $K=L=O(1/\a)$, we find a set  $X\subseteq B_\d-t$ such that
\begin{equation}\label{tmp:18.02.2011_1}
    |X| \gs \exp (-O(\epsilon^{-2} h^2 \log^2 K)) |B_{\d}| \,,
\end{equation}
and
\begin{equation}\label{tmp:18.02.2011_-1}
    | \mu_{-A} * (A+S) * \mu_{-S}(x) - 1| \ls \epsilon/3 \quad \mbox{ for all } \quad x \in hX \,.
\end{equation}
We may assume that $B_\d$ is regular.

Let $\epsilon$ be a small positive constant to be specify later.
Put $h = \lceil \log ( K/\epsilon) \rceil$ and
$l=O(\epsilon^{-4} h^2 \log^2 K)$. Applying Lemma
\ref{l:Chang+large_Bohr_exp_sums} for $X+t\sbeq B_\d$ with
parameters $\nu = O(\epsilon / (l K^{1/2}))$, $\rho = O( \epsilon /
(l d^2 K^{1/2}))$, we obtain
\begin{equation}\label{tmp:18.02.2011_2}
    |1-\gamma(x)| \ls \epsilon /(3 K^{1/2}) \quad \mbox{ for all } \quad x\in B_{\d \rho} \wedge B'_\nu
        \quad \mbox{ and } \quad \gamma \in \Spec_\epsilon (\mu_X) \,.
\end{equation}
We have $\dim (B_{\d \rho} \wedge B'_\nu)= d +O( \log^4
(1/\a))$.

By the same argument, applied for sets $A',S'$ there are sets
$X'$,
 $\L'$ of cardinality $l$ and a Bohr set $B^*_{\nu}$
that satisfy inequalities (\ref{tmp:18.02.2011_1}) and
(\ref{tmp:18.02.2011_2}), respectively. Finally, we set
$$B''= B_{\d\rho}\wedge B'_\nu \wedge B^*_{\nu}.$$
Clearly,
$d''=\dim B''=d+O( \log^4 (1/\a))$ and by Lemma
\ref{l:Bohr_est}, Lemma \ref{l:entropy_Bohr}, Lemma \ref{l:Bohr_intersection_Sanders}  and $\epsilon=\Omega ( 1)$ we have
\begin{equation}\label{in:bohr size}
 |\t{B}| \gs \exp (-O(d\log d + d\log (1/\eps) + \log^4 (1/\a)\log d +\log^5 (1/\a)+d\log(1/\a))) |B|.
\end{equation}

In view of the inequality
$$\sum_{\gamma} |\F{(A+S)}(\gamma)\F{\mu}_A(\gamma)\F{\mu}_S(\gamma)|\ls
\frac{(|A+S||A|)^{1/2}}{|S|}\ls K^{1/2},$$ which follows from
Cauchy-Schwarz inequality and Parseval's formula, we may proceed in
the same way as in the proof of Lemma 9.2 in
\cite{Sanders_2A-2A_new_optimal} and conclude that for any
probability measure $\mu$ supported on $B''$ we have
\begin{equation}\label{tmp:18.02.2011_3}
    \| (A+S) * \mu \|_{\infty} \gs 1-\epsilon \quad \mbox{ and } \quad \| (A'+S') * \mu \|_{\infty} \gs 1-\epsilon \,.
\end{equation}
Let $\eta = 1/4 d''$. We show that $(A-A')+(S-S')$ contains a
translation of $\t{B} := B''_\eta$.

 Indeed, note that

$$
    B''_{1/2} \subseteq B''_{1/2+\eta} \subseteq \dots \subseteq B''_{1/2+2d''\eta} = B'' \,.
$$
so that by pigeonhole principle, there is some $i\le 2d''$ such
that $|B''_{1/2+i\eta}| \ls \sqrt{2} |B''_{1/2+(i-1)\eta}|$. We
apply (\ref{tmp:18.02.2011_3}) for
$$
    \mu = \frac{B''_{1/2+i\eta}+B''_{1/2+(i-1)\eta}}{|B''_{1/2+i\eta}| + |B''_{1/2+(i-1)\eta}|} \,.
$$
Thus, there is $x$ such that
$$
    |(x+A+S) \cap B''_{1/2+i\eta}| + |(x+A+S) \cap B''_{1/2+(i-1)\eta}|
        \gs
            (1-\epsilon) \l( |B''_{1/2+i\eta}| + |B''_{1/2+(i-1)\eta}| \r) \,.
$$
Taking $\epsilon$ sufficiently small (see \cite{Sanders_2A-2A_new_optimal} for details), we get
$$
    |(x+A+S) \cap B''_{1/2+i\eta}| \gs \frac{3}{4} |B''_{1/2+(i-1)\eta}| \,, \quad
    |(x+A+S) \cap B''_{1/2+(i-1)\eta}| \gs \frac{3}{4} |B''_{1/2+(i-1)\eta}| \,.
$$
Analogously, for some $y$, we obtain
$$
    |(y+A'+S') \cap B''_{1/2+i\eta}| \gs \frac{3}{4} |B''_{1/2+(i-1)\eta}| \,, \quad
    |(y+A'+S') \cap B''_{1/2+(i-1)\eta}| \gs \frac{3}{4} |B''_{1/2+(i-1)\eta}| \,.
$$
Hence for each $b\in \t{B}$, we have
\begin{eqnarray*}
    (A+S) * (-A'-S') (b+y-x) &=& (x+A+S) * (-y-A'-S') (b)\\
&\gs&
            ((x+A+S) \cap B''_{1/2+i\eta}) * ((-y-A-S) \cap B''_{1/2+(i-1)\eta})
            (b)\\
       & \gs&
            |(x+A+S) \cap B''_{1/2+i\eta}| + |(y+A+S) \cap
            B''_{1/2+(i-1)\eta}|\\
                    &-&
                |((x+A+S) \cap B''_{1/2+i\eta}) \cap ((-y-A-S) \cap
                B''_{1/2+i\eta})|\\
                    &\gs&
                        \frac{3}{2} |B''_{1/2+(i-1)\eta}| - |B''_{1/2+i\eta}|
                            > 0 \,.
\end{eqnarray*}
Therefore, $(A-A')+(S-S')$ contains a translation of $\t{B}$.
Finally, by Lemma \ref{l:Reg_B}, there is $1/2\ls \sigma \ls 1$
such that $\t B_{\sigma}$ is regular. By (\ref{in:bohr size}) and Lemma \ref{l:entropy_Bohr} $\t B_\sigma$ also satisfies
(\ref{f:t(B)_size}). This completes the proof.
$\hfill\Box$

\section{Proof of the main result}
\label{sec:proof}

Let $A\sbeq \{1,\dots,N\}$ be a set having no solution to
(\ref{e}). As usually, we embed $A$ in $\Z_p$ with $p$ between
$(\sum|a_i|)N$ and $2(\sum|a_i|)N$, so $A$ has no solution to
(\ref{e}) in $\Z_p$. All sets considered below  are subsets of
$\Z_p.$ We start with the following simple observation.

\bigskip

\Lemma
{\it
 Let $B$ be a regular Bohr set of dimension $d$,   $B' \le B_\rho$ be a Bohr set and
 $\rho\ls \a/(1600d)$. Suppose that  $\mu_B (A),  \mu_B (A')\gs
\a$. Then there exists $x\in B$ such that
 \begin{equation}\label{f:l_A_and_-A_pm}
(\mu_{B'} * A) (x),\, (\mu_{B'} * A') (-x) \gs \a/4
 \end{equation}
 or
 \begin{equation}\label{f:l_A_and_-A_inc}
 \| \mu_{B'} * A \|_\infty \gs 1.5 \a  \text{~~or~~}  \| \mu_{B'} * A' \|_\infty \gs 1.5 \a
 \end{equation}
}
\label{l:A_and_-A}
\Proof By regularity of $B$ we have
$$
\a \ls \sum_{x\in B} \mu_{B} (x)  A(x) \ls \a/8 + \sum_{x\in B} (\mu_{B} * \mu_{{B'}}) (x)  A(x)
\ls \a/8+
            \frac{1}{|B|} \sum_{x\in B} (\mu_{{B'}} *  A) (x)$$
and
$$\a \ls  \a/8+\frac{1}{|B|} \sum_{x\in B} (\mu_{{B'}} *  A') (x).
$$
Hence
$$
 \sum_{x\in B} \l( (\mu_{B'} * A) (x) + (\mu_{B'} * A') (-x) \r) \gs (7\a/4)|B|
$$
and the result follows. $\hfill\Box$

\bigskip

Theorem \ref{t:Roth_Schoen_Gr''_intr} is a
 consequence of the next lemma.

\bigskip

\Lemma {\it
    Suppose that $B$ is a
    regular Bohr set of dimension $d$ and $A\sbeq B,\,$
    $\mu_B (A) = \a$ has no solution with distinct
    elements to ({\ref{e}}). Assume that
    \begin{equation}\label{l:large bohr}
    |B|\gs \exp(O(d\log d+\log^5(1/\a)+d\log(1/\a) +\log d \log^4 (1/\a)+d\log k)).
    \end{equation}
    Then there exists a regular Bohr set $B'$, such that
    \begin{equation}\label{f:Roth_Schoen_Gr_density'}
        \| \mu_{B'} * A \|_\infty \gs (1+1/(16k)) \a \,,
    \end{equation}
    $\dim B' = d + O(\log^4 (1/\a))$,
    and
    \begin{equation}\label{f:Roth_Schoen_Gr_eps'}
        |B'| \gs \exp (-O(d\log d+ \log^5 (1/\a)+d\log(1/\a) +\log d \log^4 (1/\a))) |B| \,.
    \end{equation}
}
\label{l:Roth_Schoen_Gr'}
\proof We start with mimicking the argument used by Sanders in
\cite{Sanders_Roth_1-eps}. Suppose that $\eps=c\a/(100Mdk^2)$, where
$c>0$ is a small constant and $M=\prod |a_i|$, is such that
$B_{\eps}$ is a regular Bohr set and put $B^i=(\prod_{j\not=i}
a_j)\cdot B_\eps$. By Lemma \ref{l:L_pm}, we have
\begin{equation}\label{f:multiplication_trick}
    \| k \cdot (A* \mu_B) -\sum_{i=1}^k A * \mu_B * \mu_{B^i}\|_\infty \ls 2kc\a \,.
\end{equation}
Thus, for  $\e=1/(16k),\,$  either we have  $\| \mu_{B^i} * A
\|_\infty \gs (1+\e) \a$  for some $1\ls i \ls k,$ or there is $w\in
B$ such that $\mu_{B^i}(A+w)=\mu_{B'}(a_i\cdot (A+w)) \gs (1-k\e)\a$
for every $i,$ where $B'=(\prod a_j)\cdot B_\eps$.
 In the first case we are done, so assume that
the last inequalities hold.
Since (\ref{e}) is an invariant equation we may translate our set and
assume that $\mu_{B'} (a_i\cdot A) \gs (1-k\e)\a$ for all $1\ls i\ls k$.

Let $B'_{\eps/2}\sbeq B'' \subseteq B'_\eps$ and $B''_{\eps/2}\sbeq
B''' \subseteq B''_\eps$ are regular Bohr sets. By regularity of
$B'$ and Lemma \ref{l:A_and_-A} either (\ref{f:l_A_and_-A_inc})
holds, and we are done, or there exists $x\in B'$ with
$\mu_{B''+x}(a_1\cdot A)\gs \a/8$ and $\mu_{B''-x}(a_2\cdot A)\gs
\a/8.$ We show  that there are disjoint sets $A_1,A_2$ of $A$ such
that $\a/32\ls \mu_{B''+x}(a_1\cdot A_1)\ls \a/16$ and $\a/32\ls
\mu_{B''-x}(a_2\cdot A_2)\ls \a/16$. Indeed, let $Q_1 = \{ q\in A
~:~ a_1 \cdot q \in B''+x\}$, $Q_2 = \{ q\in A ~:~ a_2 \cdot q \in
B''-x\}$. Note that $|Q_1|, |Q_2| \ge \a|B''|/8$. If $|Q_1 \cap Q_2|
> \a|B''|/16$ then split $Q_1 \cap Q_2$ into two parts $A_1$, $A_2$
whose sizes differ by at most one. Otherwise, we  put $A_1 = Q_1
\setminus Q_2$ and $A_2 = Q_2 \setminus Q_1$.

Put
$A'=A\setminus (A_1\cup A_2)$ then  $\mu_{B'} (a_i\cdot A') \gs 3\a/4$ for $i \ge 3$.
Again applying Lemma \ref{l:A_and_-A} for $B'''$
and  the arguments above,
we find $y\in B'$ and disjoint sets  $A_3,A_4\sbeq A'$ such that
$\mu_{B'''+y}(a_3\cdot A_3)\gs \a/16$ and
$\mu_{B'''-y}(a_4\cdot A_4)\gs \a/16.$

Assume that $k$ is even.
Let $l = (k-6)/2 \ge 0$. Using the arguments as before, we infer
that then there are disjoint sets $A_5,\dots,A_{k-2}$ and elements
$y_1,\dots,y_l$ such that
$$
    a_5 \cdot A_5 - y_1 \sbeq B'',~ -a_6 \cdot A_6 - y_1 \sbeq B'',~
        \dots ,~ a_{k-3} \cdot A_{k-3} - y_{l} \sbeq B'',~ -a_{k-2} \cdot A_{k-2} - y_l \sbeq B''
$$
and
$$
    \mu_{B''+y_1} (a_5 \cdot A_5),\, \mu_{B''-y_1} (a_6 \cdot A_6)\,, \dots \,,
    \mu_{B''+y_l} (a_{k-3} \cdot A_{k-3}),\, \mu_{B''-y_l} (a_{k-2} \cdot A_{k-2}) \ge \frac{\a}{16 k} \,.
$$

Finally, by Theorem \ref{t:Sanders_2A-2A_reformulation} applied to
sets
$$a_1\cdot A_1-x\sbeq B'',~ -a_2\cdot A_2-x\sbeq B'',~ a_3\cdot A_3-y\sbeq B''',~
-a_4\cdot A_4-y\sbeq B''',$$
there exists a Bohr set $\t{B}\le B'$ and $z$ such
that
\begin{equation}\label{tmp:02.05.2011_1}
    \t{B} + z \subseteq a_1\cdot A_1 +a_2\cdot A_2 + a_3\cdot A_3 + a_4\cdot A_4  + \sum_{j=5}^{k-2} a_j \cdot A_j \,,
\end{equation}
$\t{d}=\dim \t{B} = d + O(\log^4 (1/\a))$ and
\begin{equation}\label{tmp:22.02.2011_1'}
    |\t{B}| \gs \exp (-O(d\log d+ \log^5 (1/\a)+d\log(1/\a) +\log d \log^4 (1/\a))) |B| \,.
\end{equation}
The sum over $j$ in (\ref{tmp:02.05.2011_1}) can be empty.
In the case we put the sum to be equal to zero.
Notice that $z\in 4B'' + (k-6)B'''\sbeq kB''.$
Since $A_1,\dots, A_{k-2}$ are disjoint it follows that
\begin{equation}\label{tmp:13_02_2011_1'}
    a_{k-1}x_{k-1}+a_kx_k \notin \t{B} - z
\end{equation}
 for all  distinct
$ x_{k-1},x_k \in A\setminus  \cup_{j=1}^{k-2} A_j.$

By Lemma \ref{l:Reg_B} we find $1/(400k\t{d})\ls \d\ls
1/(200k\t{d})$ such that $\t{B}_\d$ is regular. Obviously
$\t{B}_{\d}$ satisfies (\ref{tmp:22.02.2011_1'}). Write
$$
    E_i:= \{ x \in B' ~:~ (\mu_{\t{B}_{\d}} * (a_i\cdot A)) (x) \gs  k/|\t{B}_{\d}| \} \,.
$$
Observe that if $-z\in E_{k-1}+E_k,$ then one  can find a solution
to (\ref{e}) with distinct $ x_1,\dots,x_k \in A.$ Therefore
$E_{k-1}\sbeq B'\setminus (-E_k-z),$ so that
$$|E_{k-1}|\ls |B'\setminus (-E_k-z)| = |B'| - |B'\cap (E_k+z)|
\ls  |B'| - |E_k| + 100 \eps M dk|B'|  \,.$$
Finally
$$|E_{k-1}|+|E_k|\ls (3/2)|B'|,$$
so that $|E_i|\ls (3/4)|{B'}|$ for some $i$. Thus
$$|A|=\|\mu_{\t{B}_{\d}} * (a_i\cdot A)\|_1\ls \|\mu_{\t{B}_{\d}} * (a_i\cdot A)\|_\infty |E_i|+(k/|\t{B}_{\d}|)|B'_{1+\d}|\,.$$
By (\ref{l:large bohr})
$$|{\t{B}_{\d}}|\gs \exp(O(-(d\log d+\log^5(1/\a)+d\log(1/\a) +\log d \log^4 (1/\a))))|B|\gs 10\cdot 8^{d+1}k/\a\,,$$
and since Lemma \ref{l:entropy_Bohr} implies  $|B'_{1+\d}|\ls
2|B'|,$ so that
$$\|\mu_{\t{B}_{\d}} * (a_i\cdot A)\|_\infty|E_i||\t{B}_{\d}|\gs 0.9 |\t{B}_{\d}||A|.$$
Hence
$$\|\mu_{B^*} *  A\|_\infty\gs 1.1\a,$$
where $B^*=a_i^{-1}\cdot {\t{B}_{\d}},$
and the assertion follows.

Now suppose that $k$ is odd. Only the first part of the proof needs to be slightly modified.
Certainly, we may assume that $a_5 = 1$.
By regularity of $B$ we have
\begin{equation}\label{f:multiplication_trick'}
    \| k \cdot (A* \mu_B) - A * \mu_B * \mu_{B''} - \sum_{i\not=5}^{} A * \mu_B * \mu_{B^i}\|_\infty \ls 2kc\a \,,
\end{equation}
where $B^i$ and  $B''$ are defined  as before. Put $l=(k-7)/2 \ge
0$. By Lemma \ref{l:A_and_-A} there are  disjoint sets $A_1,\dots,
A_k$ and elements $x,y,y_1,\dots,y_l$ such that
(\ref{tmp:02.05.2011_1})--(\ref{tmp:13_02_2011_1'}) hold. However,
$A_5\sbeq B'',$ so that $z\in kB''$. One can finish the proof in
exactly the same way as before. $\hfill\Box$

\bigskip

{\it Proof of Theorem \ref{t:Roth_Schoen_Gr''_intr}} Let $A\sbeq
B^0=\Z_p, |A|\gs \a p.$ We apply iteratively Lemma
\ref{l:Roth_Schoen_Gr'}. After $t$ steps we obtain a regular Bohr
set $B^t$ and $x_t\in \Z_p$ such that  $|A\cap (B^t+x_t)|\gs
(1+1/(16k))^t\a |B^t|,\, \dim B^t\ll t\log^4(1/\a),$ and
$$|B^t|\gs \exp(-O(t\log^4 (1/\a)\log \log (1/\a)+\log^5(1/\a)))|B^{t-1}|.$$
Since the density is always
less than $1$ we may apply Lemma \ref{l:Roth_Schoen_Gr'} at most
$O(\log(1/\a))$ times. Therefore, after  $t=O(\log (1/\a))$ iterates
assumption of Lemma \ref{l:Roth_Schoen_Gr'} are violated, so that
$$\exp(-O(\log^6 (1/\a)\log \log (1/\a)))p
\ls |B^t|\ls  \exp(O(\log^5(1/\a))),$$
which yields
$$\a\ll \exp(-c(\log p/\log \log p)^{1/6}),$$
and the assertion follows.
$\hfill\Box$

\section{Polynomial Freiman--Ruzsa Conjecture and linear equation}
\label{sec:PFRC}

Freiman-Ruzsa Polynomial Conjecture can be formulated in the following way.

\bigskip

\con\label{con1}
 {\it Let  $A\sbeq \z_N, |A|=\a N,$ then there exists a Bohr
set $B(\G,\eps)\sbeq 2A-2A$ such that $|\G|=d\ll \log(1/\a)$ and
$\eps\gg 1/\log (1/\a).$}

\bigskip

We have
$$|B(\G,\eps)|\gs \frac12\eps^{d}N,$$
so that it would give a nontrivial result provided that $\a\gg
N^{-c/\log\log N}.$ However,  it was proved in \cite{Shkredov1} and
\cite{Shkredov2} that in Chang's lemma (see section \ref{sec:Sanders_Bohr})
one can take much larger
$\eps.$ This give a (little) support for the following  version of
the above conjecture  for sparse sets.

\bigskip

\con\label{con2} {\it Let $A, A'\sbeq \z_N, |A|, |A'|\gs N^{1-c},$ then there
exists a $\d_c \log N-$dimensional Bohr set $B\sbeq A-A +A'-A'$ such that $|B|\gg N^{1-c'}$ and
$\d_c \rightarrow 0,c'\rightarrow 0$ with
$c\rightarrow 0.$ Furthermore, each $b\in B$ has
$\gg |A|^2|A'|^2/N$ representations in the form $a-b+a'-b',\, a,b\in A,\, a',b'\in A'.$ }

\bigskip

We shall give here an application of Conjecture \ref{con2}. First we recall
some definitions from \cite{Ruzsa_equations}. Let
\begin{equation}\label{equation}
a_1x_1+\dots+a_kx_k=0
\end{equation}
be an invariant linear equation.  We say that the solution $x_1,\dots,x_k$
of (\ref{equation}) is trivial if there is a partition
$\{1,\dots,k\}=\T_1\cup\dots \cup \T_{l}$ into nonempty and disjoint
sets $\T _j$ such that $x_u=x_v$ if and only if $u,v\in \T_j$ for
some $j$ and
$$\sum_{i\in \T_j}a_i=0,$$
for every $1\ls j\ls l.$ The {\it genus} of (\ref{equation}) is the
largest ${\s}$ such that there is a partition
$\{1,\dots,k\}=\T _1\cup\dots \cup \T_{\s}$ into nonempty and
disjoint sets $\T _j$ such that
$$\sum_{i\in \T_j}a_i=0,$$
for every $1\ls j\ls \s.$ Let $r(N)$ be the maximum size of a set
$A\sbeq \{1,\dots ,N\}$ having no nontrivial solution to
(\ref{equation}) with $x_i\in A$ and let $R(N)$ be the
analogous maximum over sets that the equation (\ref{equation}) has
no solution with distinct $x_i\in A$. It is not hard to prove that
$r(N)\ll N^{1/\s}.$ Much less is known about the behavior of
$R(N).$ Bukh \cite{Bukh} showed that we always have $R(N)\ll
N^{1/2-\eps}$ for the symmetric equations
$$a_1x_1+\dots+a_lx_l=a_1y_1+\dots+a_ly_l.$$
Our result is the following.

\Th {\it Assuming Conjecture \ref{con2} we have
$$R(N)\ll N^{1-c},$$
for every invariant equation (\ref{equation})
with $a_1=-a_2,~a_3=-a_4,$ where $c=c(a_1,\dots,a_k).$ }

\proof Suppose that $A$ has no solution to an equation
(\ref{equation}) with $a_1=-a_2,~a_3=-a_4,$ where
$c=c(a_1,\dots,a_k)$ and assume that $|A|\gg N^{1-c},\, c>0.$ We
embed $A$ in $\Z_M$ with $M=SN,$ where $S=\sum|a_i|,$ so that any
solution to (\ref{equation}) in $\Z_M$ is a genuine solution in
$\Z.$ Let $A=A_1\cup A_2$ be a partition of $A$ into roughly equal
parts. If Conjecture \ref{con2} holds, then there is a Bohr set
$$B\sbeq a_1\cdot A_1-a_1\cdot A_1+a_3\cdot A_1-a_3\cdot A_1$$
of dimension at most $\d_c \log N$ and  size at least $\gg N^{1-c'}.$ Put $B'=B_{1/S}.$
We show that for every $t\in \Z_M$ we have
$$|(t+B')\cap A_2|\ls k-4.$$ Indeed, if there are distinct $x_5,\dots,x_k\in (t+B')\cap A_2$, then
 $$\sum_{i=5}^k a_ix_i\in \Big (\sum_{i=5}^k a_it\Big)+B=B.$$
However, each element in $B$ has at least $|A|^4/M$ representations
in the form $a_1x-a_1y+a_3z-a_3w,\, x,y,z,w\in A_1$. This would give
a solution to (\ref{equation}) with distinct integers. Hence,
$$|B'||A_2|=\sum_t |(t+B')\cap A_2|\ls kM$$
so
$$|A|\ls 2kSN/|B'|.$$
Now, by Lemma \ref{l:entropy_Bohr} it follows that $|B'|\gg
S^{-4d}|B|\gg N^{1-c'-2\d_c \log S}.$  This leads to a contradiction,
provided  $c$ is small enough. $ \hfill\Box$
\bigskip

{\bf Acknowledgement } We wish to thank Tom Sanders for stimulating
discussions.

\bigskip

\no{Faculty of Mathematics and Computer Science,\\ Adam Mickiewicz
University,\\ Umul\-towska 87, 61-614 Pozna\'n, Poland\\} {\tt
schoen@amu.edu.pl}
\bigskip

\no{Division of Algebra and Number Theory,\\ Steklov Mathematical
Institute,\\
ul. Gubkina, 8, Moscow, Russia, 119991\\}
{\tt ilya.shkredov@gmail.com}

\end{document}